\documentclass{article}
\usepackage[T1]{fontenc}
\usepackage[utf8]{inputenc}
\usepackage{lmodern} 
\usepackage{amsmath,amsfonts,amssymb,amsthm}
\usepackage{graphicx}
\usepackage{amsrefs}

\usepackage[a4paper,top=2.5cm,bottom=2.5cm,left=2.7cm,right=2.7cm,marginparwidth=2cm]{geometry}
\usepackage[colorlinks=true, allcolors=blue]{hyperref}
\usepackage{caption}
\usepackage{subcaption}
\usepackage{booktabs}
\usepackage[table]{xcolor}

\usepackage{xspace}

\usepackage{tikz}

\usepackage{threeparttable}
\usepackage{multirow}
\usepackage{siunitx}
\usepackage{todonotes}

\usepackage{authblk}

\DeclareMathOperator{\thepace}{pace}

\newcommand{\LBpace}[1][]{\thepace^{LB}_{#1}}
\DeclareMathOperator{\npace}{npace}

\newcommand{\NLBpace}[1][]{\npace^{LB}_{#1}}
\definecolor{darkviolet}{rgb}{0.58, 0.0, 0.83}
\newcommand{\bgr}[1]{#1}

\newcommand{\btbound}[1][]{Bounds$^{\text{BT}}_{#1}$\xspace}
\newcommand{\bttime}[1][]{Time$^{\text{BT}}_{#1}$\xspace}

\newcommand*{\fullref}[1]{\hyperref[{#1}]{Section \ref*{#1}: \nameref*{#1}}}
\newcommand{\julio}[1]{#1}

\title{Impact of domain reduction techniques in polynomial optimization: A computational study}
\author[1,2]{Ignacio Gómez-Casares\thanks{Corresponding Author: ignaciogomez.casares@usc.es}}
\author[1]{Brais González-Rodríguez}
\author[1,2]{Julio González-Díaz}
\author[2]{Pablo Rodríguez-Fernández}
\affil[1]{Department of Statistics, Mathematical Analysis and Optimization and MODESTYA Research Group, University of Santiago de Compostela, 15782 Santiago de Compostela, Spain}
\affil[2]{CITMAga (Galician Center for Mathematical Research and Technology), 15782 Santiago de Compostela, Spain}

\date{\today}

\begin{document}

\maketitle

\begin{abstract}
Domain reduction techniques are at the core of any global optimization solver for NLP or MINLP problems. In this paper, we delve into several of these techniques and assess the impact they may have in the performance of an RLT-based algorithm for polynomial optimization problems. These techniques include i)~the use of (nonlinear) conic relaxations for optimality-based bound tightening, ii)~the use of Lagrangian dual information to enhance feasibility-based bound tightening, and iii)~different strategies for branching point selection. \julio{One of this paper's main contributions is providing insights into the relative impact of these techniques with respect to each other, which we hope will guide the efforts to develop and implement this type of enhancements in other solvers.} We also explore how a solver equipped with these domain reduction enhancements can further improve its performance by using machine learning to better choose the best domain reduction approach to use on a given instance.
\end{abstract}


\smallskip


\section{Introduction}
The design and implementation of global optimization algorithms for general MINLP problems is a very active field of research, and the number of available solvers has been steadily increasing over the past years. Recently, companies behind state-of-the-art MILP solvers, such as Xpress \cite{xpress} and Gurobi \cite{gurobi}, have announced the release of new versions capable of solving general MINLP problems to certified global optimality. Convexifications and domain reduction techniques are probably the two most important elements behind the efficiency of the spatial branching required to handle the nonconvexities as can be seen, for instance, in \cite{Ryoo:1996}, \cite{Tawarmalani2004}, and \cite{Belotti2009}. \julio{The literature is full of different enhancements and variations of the core ideas in these papers and, yet, we are aware of no review or consensus on their efficacy beyond the papers in which they were introduced.} The \julio{main contribution} of this paper is precisely the joint assessment of the individual impact of different aspects of domain reduction on the performance of a global optimization algorithm. We hope our analysis can guide the efforts in the development and implementation of this type of enhancements in present and future global optimization solvers for nonconvex problems.

Two essential elements present in any spatial branching algorithm that, in a broad sense, are within the scope of domain reduction, are the selection of the branching variable and the selection of the branching point: both of them allow to control how the feasible domain is reduced at a given node and, importantly, they are applied at each and every node of the branch-and-bound tree. The choice of the branching variable is relatively well understood, partially because of the thorough research for MILP problems (refer, for instance, to the seminal works \cite{Linderoth1999} and \cite{Achterberg2005}) but also due to some specialized papers for spatial branching such as \cite{Tawarmalani2004} and \cite{Belotti2009}. In the specific context of polynomial optimization, the impact of the criterion for variable selection is thoroughly discussed in \cite{Gonzalez-Rodriguez:2023} and \cite{Ghaddar:2023}. On the other hand, the impact of the selection of the branching point has received much less attention and, although different approaches are mentioned in \cite{Ryoo:1996}, \cite{Tawarmalani2004}, and \cite{Belotti2009}, we are not aware of any numerical analysis assessing how relevant these decisions may be for the performance of the resulting algorithm. In order to bridge this gap, one of the domain reduction aspects we cover in this paper is precisely the selection of the branching point.

Despite the relevance of the selection of both the branching variable and the branching point, probably the single most important aspect in domain reduction in nonlinear programming, at least in terms of the research it has generated, is bound tightening \cites{Ryoo:1996,Tawarmalani2004,Belotti2009,Puranik:2017aa}. Bound tightening techniques are a family of methods designed to reduce the search space by adjusting the bounds of the variables. These methods are often divided into i)~optimality-based bound tightening, OBBT, in which bounds are tightened by solving a series of relaxations of minor variations of the original problem and ii)~feasibility-based bound tightening, FBBT, in which tighter bounds are deduced directly by exploiting the relations between the problem constraints and the variable bounds. 
In \cite{Puranik:2017aa}, the authors present a thorough computational analysis on the huge impact that domain reduction techniques, and bound tightening in particular, have in three state-of-the-art global optimizers for NLP and MINLP problems. \cite{Gonzalez-Rodriguez:2023} reports a similar impact in the context of an RLT-based algorithm for continuous polynomial optimization problems.

Although integrating a basic bound tightening scheme in a global optimization solver is not too demanding, there is a wide range of sophisticated approaches that have been discussed in the literature and whose implementation may be significantly more complex. \cite{Belotti2013} proposes to enhance FBBT by using convex combinations of constraints instead of only the original ones, while \cite{Ryoo:1996} and \cite{Gleixner:2017} present FBBT enhancements that use Lagrangian dual information from subproblems previously solved by the algorithm. The latter two approaches can be seen as particular cases \julio{of the duality-based domain reduction framework} presented in \cite{Tawarmalani2004}. These enhanced methods often entail an additional computational overhead node by node, which should then be overcome by the reduction in the size of the resulting branch-and-bound tree. A second contribution of this paper is devoted to analyze this trade-off. We do so by studying the impact of these sophisticated approaches when integrated into a solver that already has a basic bound tightening scheme in place. More precisely, we study, in the context of the RLT-technique for polynomial optimization \cite{Sherali:1992aa}, the impact of some of the Lagrangian-based enhancements in \cite{Ryoo:1996} and \cite{Gleixner:2017}.

A novel contribution of this paper is to study the use of (nonlinear) conic relaxations, in particular second-order cone programming (SOCP) and semidefinite programming (SDP) ones, to improve the tightening obtained by OBBT. This approach is motivated by the recent analysis in \cite{Gonzalez-Rodriguez:2022}, where \julio{the authors follow the recent trend of using conic relaxations in the context of polynomial optimization \cites{Jeyakumar2017,bonami2019,elloumi2019,burer2020}}, and show the potential of using (nonlinear) conic constraints to tighten the (linear) relaxations solved by an RLT-based algorithm. We adapt their ideas with the goal of improving the performance of the underlying bound tightening techniques. \julio{The above use of SDP relaxations to solve polynomial optimization problems differs greatly from the approach initiated by J.~B.~Lasserre in \cite{lasserre2001}, with the introduction of semidefinite relaxations corresponding to liftings of the polynomial programs into higher dimensions. It has been shown that, as one moves up in the SDP hierarchy, one eventually solves the original polynomial optimization problem. Unfortunately, this approach is currently computationally prohibitive in general, since the sizes of the SDP relaxations quickly become unmanageable (refer to \cite{nie2023} for a recent monograph on the topic).}

It is worth noting that domain reduction techniques are not only relevant for theoretical purposes or for the design of general purpose solvers. There have also been important efforts to adapt these techniques to particular classes of problems, directly connected with practical applications. Particularly important are those related to optimization of power networks in general \cite{coffrin2015} and, more specifically, to the alternating current optimal power flow problem or ACOPF, a very active field of research. \cite{Chen:2016} shows the effectiveness of domain reduction techniques for ACOPF problems, \cite{Sundar2018} presents an optimality-based bound tightening scheme to improve the bounds in different parts of the network, using a (convex) quadratic relaxation of the problem, and \cite{Shchetinin2019} proposes three new methods for reducing the domain, focusing on being efficient on large-scale grids.

Last, but not least, we also build upon the recent literature on ``learning for optimization'' \cites{Lodi2017, Bengio2021,kannan2022} and, in particular, upon the approach developed in \cite{Ghaddar:2023}, to study to what extent a machine learning framework may help to increase the potential of the different domain reduction techniques under study. We analyze, individually for each domain reduction enhancement and also jointly for all of them together, to what degree one might learn how to predict the best domain reduction approach to apply on a new instance and further improve performance.

Summarizing, the study of domain reduction techniques has been a core topic in global optimization in recent years. In this work, we focus on the impact of some of these techniques in the specific context of polynomial optimization problems, building upon an RLT-based algorithm for the analysis \cite{Sherali:1992aa}. The contributions of this paper are the following: 
\begin{enumerate}
	\item To study how sensitive the performance of a branch-and-bound algorithm may be with respect to different strategies for choosing the branching point.
	\item To study the impact on performance of some advanced enhancements of FBBT and OBBT on an algorithm already equipped with basic FBBT and OBBT functionalities. In particular, a novel contribution is the study of an OBBT implementation based on the solution of SOCP and SDP relaxations.
	\item We embed the above domain reduction enhancements in a learning framework and study its potential to deliver additional improvements on performance.
\end{enumerate}

The remainder of this paper is structured as follows: In Section~\ref{sec:preliminaries}, we provide a concise overview of the reformulation-linearization technique (RLT), optimality-based bound tightening (OBBT), feasibility-based bound tightening (FBBT), and detail our computational setup for the experiments. In Section~\ref{sec:domred}, we delve into the enhancements related to domain reduction techniques. Section~\ref{sec:ml} is dedicated to the application of machine learning to select the most effective enhancements for a given problem. We then draw our study to a close in Section~\ref{sec:conclusions}.

\section{Preliminaries}\label{sec:preliminaries}

\subsection{Polynomial optimization, RLT and bound tightening}

The \textit{reformulation-linearization technique} (RLT) for solving continuous polynomial optimization problems to global optimality was originally developed in \cite{Sherali:1992aa}. More especifically, it was designed for problems of the following form:
\begin{equation}
\begin{split}
\text{minimize} & \quad \phi_0(\mathbf{x})\\
\text{subject to}  & \quad \phi_r(\mathbf{x})\geq \beta_r, \quad r=1,2,\ldots, R_1 \\
& \quad \phi_r(\mathbf{x})=\beta_r, \quad r=R_1+1,\ldots,R\\
& \quad \mathbf{x}\in\Omega \subset \mathbb{R}^{|N|}\text{,}
\end{split}
\tag{PO}
\label{eq:PO}
\end{equation}
where $N$ denotes the set of variables, each $\phi_r(\mathbf{x})$ is a polynomial of degree $\delta_r \in \mathbb{N}$, and the region $\Omega = \lbrace \mathbf{x} \in \mathbb{R}^{|N|}: 0 \leq l_j \leq x_j \leq u_j < \infty, \, \forall j \in N \rbrace \subset \mathbb{R}^{|N|}$ is a hyperrectangle containing the feasible region. Then, $\delta=\max_{r \in \{0,\ldots,R\}} \delta_r$ is the degree of the problem and $(N, \delta)$ represents all possible monomials of degree $\delta$.

RLT is based on the construction of a linear relaxation of the polynomial problem, which is then {embedded in} a branch-and-bound scheme. Problem~\eqref{eq:PO} is relaxed by replacing each monomial of degree greater than one with a corresponding RLT variable. For example, associated with a monomial of the form $x_1x_2x_4$ one would define the RLT variable $X_{124}$. More generally, RLT variables are defined as
\begin{equation}
  X_J = \prod_{j \in J}x_j, 
  \label{eq:RLTidentity}
\end{equation}
where $J$ is a multiset containing the information about the multiplicity of each variable in the underlying monomial. This linear relaxation is then solved at each node of the branch-and-bound tree. In order to get tighter relaxations and ensure convergence, additional constraints, called bound-factor constraints, are also added {(and linearized)}.\footnote{\bgr{\cite{Sherali:1992aa} proves (see Theorem 1) that, with the inclusion of the bound-factor constraints, the RLT-based algorithm converges to a global optimum of the problem. More precisely, ```either the branch-and-bound algorithm terminates finitely with the incumbent solution being optimal to \eqref{eq:PO}, or else an infinite sequence of stages is generated such that along any infinite branch of the branch-and-bound tree, any accumulation point of the $x$-variable sequence of the linear programming iterates generated at the nodes solves \eqref{eq:PO}''.}} These constraints are given, for each pair of multisets $J_1$ and $J_2$ such that \julio{the monomial associated with $J_1\cup J_2 $ belongs to $(N,\delta)$}, by
\begin{equation*}
{F_{\delta}(J_1,J_2)=\prod_{j\in J_1}{(x_j-l_j)}\prod_{j\in J_2}{(u_j-x_j)}\geq 0.}
\end{equation*}
\julio{Given that the size of $(N,\delta)$ corresponds with $\binom{n+\delta-1}{\delta}$, the number of bound-factor constraints is often prohibitively large}. In {\cite{Dalkiran:2013aa}, the authors identify a collection of monomials, which they call $J$-sets, such that convergence to a global optimum is still guaranteed if only the bound-factor constraints associated with these monomials are considered. \julio{$J$-sets are defined as those monomials of degree greater than one present in~\ref{eq:PO} which, moreover, are not included in any other monomial of~\ref{eq:PO} (multiset inclusion)}. We adopt this approach, since it can dramatically reduce the size of the resulting relaxation (at the price of getting slightly looser relaxations).}

As discussed in the introduction, domain reduction is a core element of global optimization algorithms for nonconvex optimization problems. As can be seen in \cites{Ryoo:1996,Tawarmalani2004,Belotti2009,Puranik:2017aa}, domain reduction encompasses a wide variety of techniques, which may concern different stages of a (spatial) branch-and-bound algorithm. In particular, they may involve branching decisions such as the criteria to select the branching variable or the branching point, or approaches to tighten the bounds of the variables and reduce the search space. We briefly describe below some of these approaches, since we build upon them in the subsequent sections.

For the RLT-based algorithm on which we frame our study, branching is performed as follows: At each node, once the corresponding relaxation has been solved, a score is assigned to each variable depending on the violations of the RLT identities~\eqref{eq:RLTidentity} in which it is involved. Then, the variable with the highest score is chosen for branching, splitting the interval between its lower and upper bound at its optimal value at the current node. As we will see, changes in the approach to choose the value at which branching is performed can have a substantial impact on the performance of the algorithm.

Regarding the tightening of variable bounds, although the definition of \eqref{eq:PO} requires to have bounds for the variables, they might be unnecessarily large and lead to a hyperrectangle $\Omega$ composed mainly of infeasible points. \textit{Bound tightening} (BT) refers to a series of techniques used to reduce the size of any such hyperrectangle by removing infeasible parts of it. There are two main approaches: \textit{ optimality-based bound tightening} (OBBT) and \textit{feasibility-based bound tightening} (FBBT).

OBBT techniques typically solve two optimization problems for each variable, just by substituting the objective function with the variable at hand and then minimizing and maximizing the resulting optimization problem. This procedure gives the best possible bounds for the variable given the constraints and the bounds of the other variables. The downside is that these subproblems may be as hard to solve as the original problem, so the process is often computationally prohibitive. Because of this, OBBT is normally applied to relaxations of the original problem, easier to solve, or letting it run only for a limited amount on time. Our baseline OBBT does both things: it runs on the RLT linear relaxations and does it for a limited amount of time relative to the total time available.

Even if relaxations are used for OBBT, it can still be computationally demanding, which motivates that OBBT is often used only at the root node (the approach we take) or at a few selected points of the algorithm. FBBT, on the other hand, uses constraint propagation techniques to tighten the bounds and, since no optimization problem is solved in the process, the computational cost is typically so small that FBBT is usually applied at each and every node of the branch-and-bound tree. This is precisely what we do, with an FBBT implementation that builds upon \cite{Belotti:2012aa}, where interval arithmetic and a two phase process are used to infer bounds on the variables.

\subsection{Testing environment}\label{sec:testing}

For the numerical results presented in this paper, we use the polynomial optimization solver RAPOSa \cite{Gonzalez-Rodriguez:2023}, whose core is an RLT-based algorithm. \bgr{RAPOSa is freely available at \url{https://raposa.usc.es}}. We use three sets of instances to test the performance of the various domain reduction techniques under consideration. The first one is taken from \cite{Dalkiran:2016aa} and consists of 180 instances of randomly generated polynomial optimization problems of different degrees, number of variables, and density \julio{(proportion of monomials present in the problem with respect to the total number of monomials of degree at most $\delta$)}. The second test set comes from the well-known benchmark MINLPLib \cite{Bussieck:2003aa}, a library of mixed-integer nonlinear programming instances. We have selected from MINLPLib those instances that are~\eqref{eq:PO} problems with box-constrained and continuous variables, resulting in a total of 166 instances. The third test set comes from another well-known benchmark, QPLIB \cite{Furini:2018aa}, a library of quadratic programming instances, for which we made a selection analogous to the one made for MINLPLib, resulting in a total of 63 instances. Then, we removed 6 additional instances that lead to numerical issues in some parts of the analysis, leading to a final set of 403 instances. \bgr{Instances from the three datasets can be downloaded at \url{https://raposa.usc.es/files/DS-MINLPLib-QPLIB.zip}.}



All the executions reported in this paper have been run on the supercomputer Finisterrae~III, at Galicia Supercomputing Centre (CESGA). Specifically, we used computational nodes powered with two thirty-two-core Intel Xeon Ice Lake 8352Y CPUs with 256GB of RAM and 1TB SSD. All the executions have been run taking as stopping criterion that the relative or absolute gap is below the threshold $0.001$. The time limit of each execution was set to 1 hour and OBBT was limited to 20\% of the total time, \emph{i.e.}, 12 minutes.

We use different measures to assess the performance of the different approaches. We describe them here to avoid repetition and simplify the notation:
\begin{itemize}
 \item \textbf{Solved}. Number of instances solved to certified optimality (relative or absolute gap below $0.001$).
 \item \textbf{Gap}. Geometric mean of the optimality gap obtained by each approach. We disregard instances for which i)~at least one approach did not return an optimality gap and (ii)~all the approaches solved it within the time limit.\footnote{\julio{The choice of the geometric mean for the gap is motivated by the fact that, for some instances, $LB<0$, which means that the relative gap, computed as $(UB-LB)/|UB|$, can take values arbitrarily large. Also, there is no need to rely on the shifted mean to account for values very close to zero since, for the numerical analysis, the value of the gap for all solved instances is set to $0.001$, the stopping tolerance.}}
 \item \textbf{Time}. Geometric mean of the running time of each approach. We disregard instances for which i)~every approach solved it in less than 5 seconds and (ii)~no approach solved it within the time limit.
 \item \textbf{Pace}. Geometric mean of $\LBpace$, as introduced in \cite{Ghaddar:2023}. It represents the number of seconds needed to improve the lower bound of the algorithm by one unit (the pace at which the lower bound improved). We disregard instances solved by every approach in less than 5 seconds. As thoroughly discussed in \cite{Ghaddar:2023}, the main motivation behind this performance measure is that it allows to compare the performance on all instances together, whereas Time and Gap fail to do so. Time is not informative when comparing performance between instances not solved by any approach (all of them reach the time limit). These instances can be compared using Gap which, in turn, is not informative in instances solved by all approaches (all of them close the gap), and where Time could be used. Pace, on the other hand, is informative regardless of the number of configurations that might have solved each of the different instances. 
 \item \textbf{Nodes}. Geometric mean of the number of nodes explored in the branch-and-bound tree. We only consider the instances solved by all the approaches within the time limit.
 \item \textbf{\btbound}. Mean of the average improvement of the variable bounds after applying the corresponding bound tightening approach at the root node. 
 \item \textbf{\bttime}. Mean of the bound tightening time at the root node. 
\end{itemize}

We have chosen to report geometric means in our results, with the exception of \btbound and \bttime, for which the standard mean is used. Geometric means are known to have the advantage of being less sensitive to outliers, which often makes them preferable for performance comparisons. The reason to report the standard mean for \btbound is that it measures how much the bounds are tightened with the different approaches, but it is not really a performance measure as Gap, Time, or Pace, for instance. Thus, since \btbound measures how much certain parameters of the problem are changed before the branching starts, reporting the standard averages allows to get a better picture of the impact of the bound tightening phase on the search space of the algorithm. Similarly, \bttime is not a measure of the overall performance of a given approach, but just a measure of the time spent on bound tightening at the root node. Reporting standard averages provides a more informative quantification of the impact of a given approach during this phase of the algorithm.\footnote{Note that we do not report the number of instances for which each approach did not return an optimality gap. The reason for this is that there are not significant differences across approaches and so disregarding them simplifies the tables of results.}


\section{Domain reduction in the context of RLT}\label{sec:domred}

In this section we present different enhancements for some domain reduction techniques. First, we study the impact of second-order cone programming (SOCP) and semidefinite programming (SDP) on OBBT. Second, we analyze how duality in linear programming can improve the performance of FBBT. Finally, we measure the impact of changing the rule for selecting the branching point at the different nodes of the branch-and-bound tree.

\subsection{OBBT enhancements based on conic optimization}
\label{sec:conic}

In recent literature, such as \cite{Burer2006} and \cite{Buchheim2012}, there has been an increase in the research devoted to the use of different SOCP and SDP constraints to define tighter relaxations of a given nonlinear problem, specially in the study of polynomial optimization problems. A related approach consists in adding linearizations of the SDP constraints, called linear SDP-based cuts, as in \cite{Sherali2011}, \cite{Baltean-Lugojan:2019}, and \cite{Gonzalez-Rodriguez:2023}. In \cite{Gonzalez-Rodriguez:2022} the authors use second-order cone and semidefinite constraints in order to tighten the standard RLT linear relaxations along the branch-and-bound tree. Their results show that the use of (nonlinear) conic constraints has significant potential, particularly so for some specific subclasses of problems in the test sets they consider. In this section we build upon the approach in \cite{Gonzalez-Rodriguez:2022}, using the same collection of SOCP and SDP constraints but only to tighten the linear relaxations that OBBT has to solve at the root node, so no specialized conic optimization solver is required beyond the root node. The goal is to study to what extent this enhanced OBBT can provide tighter bounds for the variables. Since solving these conic relaxations is usually more difficult than solving the original linear relaxations, there is a trade-off between time consumption and the quality of the bounds. In order to get the best out of this trade-off, when deciding which SOCP and SDP constraints to add to tighten the linear relaxations, one has to be careful to preserve whatever sparsity structure the original problem may have. To do so, in \cite{Gonzalez-Rodriguez:2022} the authors anchor the definition of the second-order cone and semidefinite constraints to the $J$-sets of~\eqref{eq:PO} as follows:
\begin{itemize}
\item \textbf{SOCP}: For each $J$-set $J$ and for each pair of variables present in $J$, $x_i\neq x_j$, we add the following second-order cone constraint to the relaxation solved at the OBBT phase:
\begin{align}\label{eq:socp}
\frac{X_{ii} + X_{jj}}{2} &\geq \Biggl\lVert \begin{pmatrix}
    X_{ij}\\
    \displaystyle\frac{X_{ii} - X_{jj}}{2}\\
    \end{pmatrix} \Biggr\rVert_2.
\end{align}
\item \textbf{SDP$^1$}: For each $J$-set, we define the vector $\omega^1=(x_j)_{j\in J}$ and we build matrix $M=(\omega^1)^T(\omega^1)$. We denote the linearization of matrix $M$ by $[M]_L$, obtained by replacing all monomials in it by the corresponding RLT variables. Then, we add, to the linear relaxation solved at the OBBT phase, the constraint $[M]_L\succcurlyeq 0$ (positive semidefiniteness).
\item \textbf{SDP$^2$}: We do the same as in SDP$^1$ but changing the vector $\omega^1$ to $\omega^2=(1,(x_j)_{j\in J})$.
\end{itemize}

It is worth noting that, in \cite{Gonzalez-Rodriguez:2022}, the authors discuss some additional approaches which, in order to reduce the number of conic constraints added to the standard relaxations and further reduce the computational effort to solve them, check what constraints are binding at the root node and then disregard all the nonbinding ones in the subsequent relaxations. Although such an approach might also be applied to our OBBT relaxations, we believe it is not advisable, since the relaxations solved by the OBBT have completely different objective functions, so there might be no correspondence between the binding constraints of the different subproblems.

\subsubsection{Numerical results}

The default configuration of RAPOSa relies on Gurobi \cite{gurobi} for the solution of the linear relaxations. For the analysis in this section, we use Gurobi and Mosek \cite{mosekmanual} for the SOCP relaxations. The SDP ones are solved with Mosek, since Gurobi cannot handle them. With the goal of assessing the potential of an OBBT whose relaxations have been tightened with (nonlinear) conic constraints, we compare the performance of the different conic approaches with a baseline configuration already equipped with standard OBBT and FBBT functionalities:
\begin{itemize}
    \item \texttt{Baseline}: We use the standard OBBT and FBBT used in~\cite{Gonzalez-Rodriguez:2023}. \bgr{At the root node, we perform OBBT before solving it. OBBT consists of minimizing and maximizing each variable over the feasible set of the root node (\emph{i.e.}, its linear relaxation). Moreover, before solving each node, we perform FBBT using interval arithmetic.}
    \item \texttt{SOCP$^G$}: We follow approach \textbf{SOCP} above and solve the relaxations with Gurobi.
    \item \texttt{SOCP$^M$}: We follow approach \textbf{SOCP} above and solve the relaxations with Mosek.
    \item \texttt{SDP$^1$}:  We follow approach \textbf{SDP$^1$} above and solve the relaxations with Mosek.
    \item \texttt{SDP$^{2}$}: We follow approach \textbf{SDP$^2$} above and solve the relaxations with Mosek.
\end{itemize}


\begin{table}[!htpb]
    \caption{Performance of the different conic OBBT approaches.}
    \label{tab:conicbt}
    \centering
    {\small \addtolength{\tabcolsep}{-2pt} 
    \begin{tabular}{lrrrrr|rr}
    \toprule
        {\scriptsize (403 instances)} & Solved {\scriptsize (403)} & Gap {\scriptsize (111)} & Time {\scriptsize  (153)} & Pace {\scriptsize (268)} & Nodes {\scriptsize (285)} & \btbound {\scriptsize (403)} & \bttime {\scriptsize (403)} \\
        \texttt{Baseline} & 286 & 0.132 & 56.95 & 6.56 & 104.43 & 0.338 & 280.99 \\
        \texttt{SOCP$^G$} & 285 & 0.133 & 97.67 & 9.98 & \cellcolor{gray!25}104.38 & \cellcolor{gray!25}0.341 & 391.92 \\
        \texttt{SOCP$^M$} & 286 & \cellcolor{gray!25} 0.128 & 65.69 & 6.99 & \cellcolor{gray!25}102.20 & \cellcolor{gray!25}0.344 & 287.92 \\
        \texttt{SDP$^1$} & 286 & 0.150 & 114.42 & 15.19 & \cellcolor{gray!25}103.62 & \cellcolor{gray!25}0.339 & 404.08 \\
        \texttt{SDP$^{2}$} & 286 & 0.159 & 132.93 & 17.89 & \cellcolor{gray!25}103.63 & \cellcolor{gray!25}0.345 & 433.57 \\
    \bottomrule
    \end{tabular}}
\end{table}

In Table~\ref{tab:conicbt} we present the results associated to the aforementioned approaches, where shaded numbers are used to highlight whenever an approach outperforms the baseline for a certain metric. At first sight, it seems that the overall performance of the new conic approaches for OBBT is worse than the one delivered by \texttt{Baseline}, specially for the SDP ones. As expected, with respect to \texttt{Baseline}, all the conic approaches reduce the number of nodes of the branch-and-bound tree and improve the tightness of the bounds (\btbound). Yet, the computational overhead required by OBBT at the root node (\bttime) is not compensated by the modest improvement on \btbound. This is is also captured by the pace of the different approaches: although \texttt{SOCP$^G$} and \texttt{SOCP$^M$} are very competitive in terms of Gap, the fact that they are more demanding computationally than \texttt{Baseline} also makes them worse in terms of Pace. \texttt{SOCP$^M$} is clearly the best performing conic approach, being fairly competitive with \texttt{Baseline}. Indeed, if we looked at standard means instead of geometric ones (not shown in the table), we would see that \texttt{SOCP$^M$} and \texttt{Baseline} are essentially tied for both time ($368.18$ vs $369.81$) and pace ($216108.52$ vs $215531.70$). The main weakness of the conic approaches, and specially the SDP ones, is that they are more ``risky'', in the sense that, since conic solvers are not yet as efficient as linear ones, a conic approach may turn out to be particularly slow in some instances, which might overshadow promising behavior on others.


The results in Table~\ref{tab:conicbt} go along the same lines of those reported in \cite{Gonzalez-Rodriguez:2022}, where these conic enhancements are not included at the OBBT stage of the algorithm but, instead, they are added to the RLT relaxation in each and every node of the branch-and-bound tree. The authors go on to show that, despite this initially discouraging aggregate behavior, there are specific instances and subclasses of problems in the tests under study in which some of the SOCP or SDP approaches significantly outperform their baseline configurations. In order to check whether or not we are in a similar situation, in Table~\ref{tab:conicbt5} we present, for each approach under study, how often it is within $5\%$ of the value of the best performing approach according to the given metric. 

\begin{table}[!htpb]
    \caption{Frequency with which each conic OBBT approach is within 5\% of the best one.}
    \label{tab:conicbt5}
    \centering
    {\small \addtolength{\tabcolsep}{-2.5pt} 
    \begin{tabular}{lrrrr|rr}
    \toprule
        {\scriptsize (403 instances)} & Gap$_{5\%}$ {\scriptsize (111)} & Time$_{5\%}$ {\scriptsize (153)} & Pace$_{5\%}$ {\scriptsize (268)} & Nodes$_{5\%}$ {\scriptsize (285)} & \btbound[5\%] {\scriptsize (403)} & \bttime[5\%] {\scriptsize(403)}\\
        \texttt{Baseline} & 100 & 120 & 227 & 208 & 243 & 361 \\
        \texttt{SOCP$^G$} & 88 & 42 & 150 & \cellcolor{gray!25}209 & \cellcolor{gray!25} 259 & 68 \\
        \texttt{SOCP$^M$} & \cellcolor{gray!25}103 & 93 & 200 & \cellcolor{gray!25}224 & \cellcolor{gray!25} 283 & 141 \\
        \texttt{SDP$^1$} & 82 & 23 & 116 & \cellcolor{gray!25}214 & \cellcolor{gray!25} 249 & 54 \\
        \texttt{SDP$^{2}$} & 80 & 10 & 102 & \cellcolor{gray!25}230 & \cellcolor{gray!25} 358 & 37 \\
    \bottomrule
    \end{tabular}}
\end{table}

The results are consistent with those in Table~\ref{tab:conicbt}, with \texttt{Baseline} and \texttt{SOCP$^M$} coming out on top. Importantly, Table~\ref{tab:conicbt5} also shows that \texttt{Baseline} is not within 5\% of the best performing approach in a significant amount of instances: 10\% for Gap, 22\% for Time, and 15\% for Pace. Thus, we could get significant performance improvements if we were able to anticipate the best approach to run on a given instance. This idea is explored in Section~\ref{sec:ml}, where we embed these approaches in a machine learning framework, and train a model with the goal of predicting what the best approach would be on an unseen instance.

\julio{Finally, it is worth noting that, in the analysis carried out in \cite{Gonzalez-Rodriguez:2022}, the authors managed to successfully identify families of instances where one specific approach was particularly effective. Unfortunately, we have not been able to find such patterns behind the results in Table~\ref{tab:conicbt} and Table~\ref{tab:conicbt5}, or any other set of results presented in this work.} 

\subsection{FBBT enhancements based on duality}
\label{sec:bound-tightening}

In this section we take the ideas developed in \cite{Ryoo:1996} and \cite{Gleixner:2017} and study their impact when added on top of the basic bound tightening, OBBT and FBBT, discussed in Section~\ref{sec:preliminaries} for the RLT technique. Both works were developed in the context of global solvers for non-convex optimization problems, where a convex relaxation is used in a branch-and-bound scheme. They discuss different approaches to improve the performance of the FBBT stages, in the sense of obtaining even tighter bounds, which we adapt below to the context of polynomial optimization.

The core idea is the same in both \cite{Ryoo:1996} and \cite{Gleixner:2017}: make use of the Lagrangian multipliers associated to previously solved subproblems to obtain valid cuts that can be added to FBBT's constraint propagation scheme. This should lead to (weakly) tighter bounds whenever FBBT is invoked but, at the same time, entail some computational overhead. Both approaches use these new cuts to extend the set of constraints on which the FBBT is run and they do so at all nodes of the branch-and-bound tree. In \cite{Gleixner:2017} they use information captured during the initial (or previous) pass of OBBT, while in \cite{Ryoo:1996} they use information from the solution of the relaxation at the current node in the tree. For the sake of completeness, we briefly describe below the two approaches and how they have been adapted to our context.

\julio{We start by outlining the approach proposed in \cite{Gleixner:2017}. Suppose that $U$ is a valid upper bound for \eqref{eq:PO}, available at the moment of applying OBBT. We can tighten the OBBT relaxations by adding constraint $c^T x\leq U$. Then,} after solving all OBBT's problems, we have, for each variable, one solution for its maximization problem and another one for its minimization problem. Suppose we have the primal-dual solution $(\tilde{x}, \tilde{\lambda}, \tilde{\mu})$ from solving the minimization problem for variable $x_k$, where $\tilde{\lambda}$ is the vector of multipliers associated with the constraints of the OBBT relaxation and multiplier $\tilde{\mu}$ corresponds \julio{with constraint $c^T x\leq U$. Finally, let $\tilde r\in \mathbb R^n$ be the vector of reduced costs associated with $\tilde{x}$. Theorem~1 in \cite{Gleixner:2017} establishes that} the following constraint is a valid cut for~\eqref{eq:PO}:
\begin{equation}
x_k \geq \sum_{j \in \lbrace 1, \dots, n\rbrace} \tilde{r}_j x_j + \tilde{\mu} U + \tilde{\lambda}^T b\text{,}
\label{eq:OBBT1}
\end{equation}
\julio{where $b$ comes from the right-hand sides of the linear relaxation used by OBBT (including the $\beta_r$ coefficients and the right-hand sides of the bound-factor constraints)}. Similarly, we can use the maximization OBBT problem to obtain \julio{cut $x_k \leq \sum_{j \in \lbrace 1, \dots, n\rbrace} \tilde{r}_j x_j + \tilde{\mu} U + \tilde{\lambda}^T b$. Interestingly, if an improved solution of~\eqref{eq:PO} is found, Theorem~1 in \cite{Gleixner:2017} also ensures that the above cuts can be tightened by updating the value of $U$.} 

Now, we describe the cuts presented in \cite{Ryoo:1996}. Consider the relaxation of \eqref{eq:PO} at a given node. \julio{Let $L$ be its optimal value and suppose that a certain constraint $g(x) \leq 0$} is active at the solution of this node.\footnote{\julio{For instance, for the original constraints of \eqref{eq:PO}, they would be written as $\beta_r-[\phi_r(\mathbf{x})]_L\leq 0$, where $[\phi_r(\mathbf{x})]_L$ represents the linearization of $\phi_r(\mathbf{x})$.}} Suppose, moreover, that this constraint has \julio{associated the dual multiplier $\tilde{\mu}$} and that $U$ is a valid upper bound for~\eqref{eq:PO}. Then, the following constraint is a valid cut for~\eqref{eq:PO}:
\begin{equation*}
g(x) \geq -(U-L)/{\tilde{\mu}}\text{.}
\end{equation*}

\julio{The above cuts also apply} to constraints of the form $x_k - x_k^U \leq 0$ and $x_k^L - x_k \leq 0$ (where $x_k^U$ and $x_k^L$ are the upper and lower bounds of variable $x_k$ \julio{at the given node}), obtaining
\begin{equation}
 x_k \geq x_k^U - (U-L)/{\tilde{\mu}} \quad \text{and} \quad x_k \leq x_k^L + (U-L)/{\tilde{\mu}}\text{.}
 \label{eq:OBBT2}
\end{equation}

If we had an equality constraint \julio{$g(x) = 0$}, we would use the sign of $\tilde{\mu}$ to determine the sense of the new cut's inequality. \julio{Again, is worth emphasizing that the values for the bounds $U$ and $L$ in the above cuts} are not static: they are updated as the algorithm progresses. Put differently, whenever any of these constraints are added to the FBBT at a given node, the current values for the bounds are used.

\subsubsection{Numerical results}

In order to assess the performance of the aforementioned approaches, we have studied \julio{a wide range of variations in how often \texttt{FBBT$^{OB}$} and \texttt{FBBT$^{NB}$} are applied along the branch-and-bound tree}. We present here the best performing ones:

\begin{itemize}
    \item \texttt{Baseline}: We use the standard OBBT and FBBT used in~\cite{Gonzalez-Rodriguez:2023}.
    \item \texttt{FBBT$^{OB}$}: The OBBT-derived cuts from \cite{Gleixner:2017} are added to FBBT at the root node and at nodes whose depth is a multiple of 50.
    \item \texttt{FBBT$^{NB}$}: The node by node cuts from \cite{Ryoo:1996} are added to FBBT at the root node and at nodes whose depth is a multiple of 10.
    \item \texttt{FBBT$^{OB+NB}$}: Both approaches are applied together.
\end{itemize}

In Section~\ref{app:FBBT} in the Appendix we report some additional results which, in particular, show that applying these enhancements at all nodes is not desirable. This is because the computational overhead required to perform FBBT's constraint propagation on the extended set of constraints is not compensated by the gains from bound reduction. \texttt{FBBT$^{OB}$} and \texttt{FBBT$^{NB}$} provide a good balance between the computational effort and the reduction of the size of the branch-and-bound tree.


\sisetup{
	round-mode = places,
	round-precision = 2,
	detect-weight=true,
	table-number-alignment = right,
	table-alignment-mode = none,
	table-text-alignment = right
	}

\begin{table}[!htpb]
    \caption{Performance of the FBBT enhanced with cuts from \cite{Gleixner:2017} and \cite{Ryoo:1996}.}
    \label{tab:zibbaronbt}
    \centering
    {\small \addtolength{\tabcolsep}{-2pt} 
    \begin{tabular}{l r S[round-precision=3] S S S | S}
    \toprule
        {\scriptsize (403 instances)} & {Solved {\scriptsize(403)}} & {Gap {\scriptsize(111)}} & {Time {\scriptsize(134)}} & {Pace {\scriptsize(248)}} & {Nodes {\scriptsize(285)}} & {\bttime {\scriptsize(403)}} \\
        {\texttt{Baseline}} & 286 & 0.132475 & 100.984600 & 8.293355 & 104.433018 & 280.991376 \\
        {\texttt{FBBT$^{OB}$}} & \cellcolor{gray!25} 287 & \cellcolor{gray!25} 0.125897 & 102.796067 & \cellcolor{gray!25} 8.269012 & \cellcolor{gray!25} 101.222215 & 300.358402 \\
        {\texttt{FBBT$^{NB}$}} & \cellcolor{gray!25} 287 & \cellcolor{gray!25} 0.127776 & \cellcolor{gray!25} 100.476846 & \cellcolor{gray!25} 8.211333 & \cellcolor{gray!25} 98.611332 & 299.392319 \\
        {\texttt{FBBT$^{OB + NB}$}} & \cellcolor{gray!25} 287 & \cellcolor{gray!25} 0.127502 & \cellcolor{gray!25} 98.320568 & \cellcolor{gray!25} 8.084024 & \cellcolor{gray!25} 94.519323 & 317.244532 \\
    \bottomrule
    \end{tabular}}
\end{table}

The results in Table~\ref{tab:zibbaronbt} show that the duality-based approaches consistently outperform \texttt{Baseline} for the different performance measures. Yet, the improvements are relatively small and, as expected, far from the ones reported in \cite{Gonzalez-Rodriguez:2023} when basic versions of OBBT and FBBT were added to an RLT implementation without any such functionality. In their numerical analysis, they obtain 70\%-90\% improvements for Gap, Time, and Nodes in MINLPLib instances. Thus, the results in Table~\ref{tab:zibbaronbt}, combined with those in \cite{Gonzalez-Rodriguez:2023}, confirm bound tightening techniques as a fundamental element of global optimization solvers but, at the same time, they show that the effort required to implement sophisticated variants of these techniques might not always be worthwhile. In particular, performance improvements such as the ones in Table~\ref{tab:zibbaronbt} could be overshadowed by potential stability and robustness issues coming from these additional functions and unforeseen interference with other current or future functionalities.

\subsection{Impact of the branching point on the branch-and-bound tree} \label{sec:branching-point}

Whenever a subproblem is solved at a node of the branch-and-bound tree, important decisions have to be made. Two crucial ones are the selection of the branching variable and of the precise value of that variable on which branching is done. As mentioned in the introduction, the choice of the branching variable has been thoroughly discussed by past literature, mainly for MILP problems (see, for instance, \cite{Linderoth1999} and \cite{Achterberg2005}), but also in the context of nonlinear optimization and spatial branching in \cite{Tawarmalani2004} and \cite{Belotti2009}.  When it comes to polynomial optimization and the RLT technique in particular,  \cites{Gonzalez-Rodriguez:2023,Ghaddar:2023,Ghaddar:2025} offer deep analyses on the impact of variable selection.\footnote{\julio{In the seminal paper \cite{Sherali:1992aa}, starting from an optimal solution $(\bar X, \bar x)$ at a node in the branch-and-bound tree, the violation of each variable is computed as $\theta_j = \max_{J\subset (N,\delta)\colon |J|<\delta} \lbrace \vert \bar{X}_{J \cup \lbrace j \rbrace} - \bar{x}_j \bar{X}_{J} \vert$, where $J \cup \{j\}$ considers only the monomials present in~\ref{eq:PO}. Then, the variable that maximizes this value is chosen for branching. More recent contributions rely on variations where the maxima in the computation of the $\theta_j$ values are replaced with weighted sums.}}

The selection of the branching point, on the other hand, has not received much attention by past literature and, importantly, we have seen no numerical analysis studying how the performance of a given solver might be affected by different strategies, which is precisely the main goal of this section. \cite{Ryoo:1996} and \cite{Belotti2009} briefly discuss the main strategies for selecting the branching point in the state-of-the-art solvers Couenne \cite{Belotti2009} and BARON \cite{Sahinidis2024}, respectively. Suppose that a branching variable has already been chosen after solving the subproblem at a certain node. Then, given the current range of that variable, \emph{i.e.}, the interval between its current lower and upper bounds, different natural branching points can be considered: the middle point of the interval, the value of the variable in the optimal solution of the subproblem, or some convex combination of these two points. \cite{Ryoo:1996} and \cite{Belotti2009} also suggest to branch on the value of the variable in the best solution already found by the algorithm (provided it belongs to the current range of the variable).

\subsubsection{Numerical results}\label{sec:branching-point-numerical}
We now present and compare a few simple strategies for the selection of the branching point, which build upon the ideas outlined in \cite{Ryoo:1996} and \cite{Belotti2009}. Given the branching variable, we denote by $OV$ its optimal value at the current node and by $MP$ the middle point of its current range. Since, in some preliminary analysis, we observed that branching whenever possible on the value of the variable in the best available solution yields slight performance improvements at no cost, this approach is incorporated in the following five approaches:

\begin{itemize}
    \item \texttt{Baseline}: Branch on $OV$.
    \item \texttt{0.75$\cdot$OV+0.25$\cdot$MP}: Branch on the convex combination $0.75\cdot OV + 0.25\cdot MP$.
    \item \texttt{0.5$\cdot$OV+0.5$\cdot$MP}: Branch on the convex combination $0.5\cdot OV + 0.5\cdot MP$.
    \item \texttt{0.25$\cdot$OV+0.75$\cdot$MP}: Branch on the convex combination $0.25\cdot OV + 0.75\cdot MP$.
    \item \texttt{MP}: Branch on $MP$.
\end{itemize}

It is worth noting that we also studied some other slightly more complex strategies along the lines suggested in \cite{Ryoo:1996} and \cite{Belotti2009}, but they did not result in extra performance improvements. Thus, we chose to stick to these simple strategies for the sake of exposition.

\begin{table}[!htpb]
    \caption{Performance of different approaches for branching point selection.}
    \label{tab:branchingpoint}
    \centering
    {\small \addtolength{\tabcolsep}{-2pt} 
    \begin{tabular}{l r *{5}{r}}
    \toprule
        {\scriptsize (403 instances)} & {Solved {\scriptsize(403)}} & Gap {\scriptsize(110)} & Time {\scriptsize(138)} & Pace {\scriptsize(250)} & Nodes {\scriptsize(280)} \\
        \texttt{Baseline} & 286 & 0.124 & 99.21 &8.16 & 105.70 \\
        \texttt{0.75$\cdot$OV+0.25$\cdot$MP} & \cellcolor{gray!25} 287 & \cellcolor{gray!25} 0.111 & \cellcolor{gray!25} 93.22 & \cellcolor{gray!25} 7.68 & \cellcolor{gray!25} 99.58 \\
        \texttt{0.5$\cdot$OV+0.5$\cdot$MP} & \cellcolor{gray!25} 288 & \cellcolor{gray!25} 0.108 & \cellcolor{gray!25} 96.75 & \cellcolor{gray!25} 7.64 & \cellcolor{gray!25} 102.10 \\
        \texttt{0.25$\cdot$OV+0.75$\cdot$MP} & \cellcolor{gray!25} 287 & \cellcolor{gray!25} 0.108 & \cellcolor{gray!25} 97.79 & \cellcolor{gray!25} 7.70 & \cellcolor{gray!25} 105.23 \\
        \texttt{MP} & \cellcolor{gray!25} 287 & \cellcolor{gray!25} 0.109 & 107.71 & \cellcolor{gray!25} 8.09 & 107.98 \\
    \bottomrule
    \end{tabular}}
\end{table}

In Table~\ref{tab:branchingpoint} we see that all the new approaches beat \texttt{Baseline}, which is the worst one according to Solved, Gap, and Pace. Further, \texttt{Baseline} is second to last according to Time and Nodes, only ahead of \texttt{MP}. \texttt{Baseline} and \texttt{MP} are significantly outperformed by the strict convex combinations, which shows that both \texttt{OV} and \texttt{MP} should be taken into account to determine the branching point. Although the three strict convex combinations deliver similar results, \texttt{0.5$\cdot$OV+0.5$\cdot$MP} seems to be slightly superior, as it comes out on top in Solved, Gap, and Pace. 

Given the results in Table~\ref{tab:branchingpoint}, we see that \texttt{OV} and \texttt{MP} are particularly bad at Gap and Time, respectively. \julio{These results suggest that sticking to the ``extreme'' approaches, \texttt{OV} or \texttt{MP}, may be detrimental in performance with respect to the ``interior'' approaches.} Put differently, would a hypothetical version of the solver capable of choosing the best performing approach between \texttt{OV} and \texttt{MP} outperform the strict convex combinations? The answer to this question is negative. Table~\ref{tab:branchingpoint5} shows, for instance, that the aforementioned hypothetical configuration would be within 5\% of the best performing approach in Time in, at most, 96 out of the 138 reported instances, whereas approach \texttt{0.5$\cdot$OV+0.5$\cdot$MP} on its own is within 5\% of the best performing one in 102 instances.

\begin{table}[!htpb]
    \caption{Frequency with which each branching point approach is within 5\% of the best one.}
    \label{tab:branchingpoint5}
    \centering
    {\small \addtolength{\tabcolsep}{-2pt} 
    \begin{tabular}{lrrrr}
    \toprule
        {\scriptsize (403 instances)} & Gap$_{5\%}$ {\scriptsize (110)} & Time$_{5\%}$ {\scriptsize (138)} & Pace$_{5\%}$ {\scriptsize (250)} & Nodes$_{5\%}$ {\scriptsize (280)} \\
        \texttt{Baseline} & 65 & 59 & 137 & 151  \\
        \texttt{0.75$\cdot$OV+0.25$\cdot$MP} & \cellcolor{gray!25}82 & \cellcolor{gray!25}102 & \cellcolor{gray!25}193 & \cellcolor{gray!25}207  \\
        \texttt{0.5$\cdot$OV+0.5$\cdot$MP} & \cellcolor{gray!25}100 & \cellcolor{gray!25}81 & \cellcolor{gray!25}180 & \cellcolor{gray!25}195  \\
        \texttt{0.25$\cdot$OV+0.75$\cdot$MP} & \cellcolor{gray!25}99 &\cellcolor{gray!25}72 & \cellcolor{gray!25}174 & \cellcolor{gray!25}158  \\
        \texttt{MP} & \cellcolor{gray!25}97 & 37 & \cellcolor{gray!25}140 & 119  \\
    \bottomrule
    \end{tabular}}
\end{table}

\subsection{Impact of combined enhancements}\label{sec:combined}
As a complement to the individual analysis in the preceding subsections, regarding the performance of the different approaches to enhance domain reduction techniques, we now study the impact of incorporating all of them together. Note that, although all the enhancements deal with different aspects of domain reduction, they are essentially independent from one another. \julio{Thus, a priori, there is no reason to expect that combining them might have a detrimental effect on performance}: i)~the enhancements in Section~\ref{sec:conic} only affect the OBBT at the root node, ii)~the enhancements in Section~\ref{sec:bound-tightening} just affect the FBBT at all nodes, and iii)~the enhancements in Section~\ref{sec:branching-point} deal with the branching point selection, which is independent of the OBBT and FBBT functionalities. For the sake of exposition, we just report the performance of combinations of the best performing approaches for each individual enhancement: \texttt{SOCP$^M$}, \texttt{FBBT$^{OB + NB}$}, and \texttt{0.5$\cdot$OV+0.5$\cdot$MP}. Yet, it is important to note that retrieving dual information for (nonlinear) conic problems is far from straightforward. Thus, whenever a combination uses \texttt{SOCP$^M$}, we replace \texttt{FBBT$^{OB + NB}$} with \texttt{FBBT$^{NB}$}, the best performing approach in Section~\ref{sec:bound-tightening} once we disregard the approaches applying \texttt{FBBT$^{OB}$}.

\sisetup{
	round-mode = places,
	round-precision = 2,
	detect-weight=true,
	table-number-alignment = right,
	table-alignment-mode = none,
	table-text-alignment = right
	}

\begin{table}[!htpb]
    \caption{Performance of RAPOSa combining the enhancements.}
    \label{tab:combined}
    \centering
    {\small \addtolength{\tabcolsep}{-2pt} 
    \begin{tabular}{l r S[round-precision=3] S S}
    \toprule
    {\scriptsize (403 instances)} & {Solved {\scriptsize (403)}} & {Gap {\scriptsize (112)}} & {Time {\scriptsize (146)}} & {Pace {\scriptsize (257)}} \\
    \midrule
	\texttt{Baseline}	                                                     & 286                       & 0.119948	                    &	83.073973		            &	7.619640	\\
 	\texttt{SOCP$^M$}	                                                     & 286	                      & \cellcolor{gray!25} 0.116285	&	93.331352		            &	8.045351	\\
	\texttt{FBBT$^{OB + NB}$}	                                             & \cellcolor{gray!25} 287 	                    &\cellcolor{gray!25} 0.115487 	  &\cellcolor{gray!25} 81.449627	&\cellcolor{gray!25} 7.454167	\\
	\texttt{0.5$\cdot$OV+0.5$\cdot$MP}	                                     & \cellcolor{gray!25} 288                       &\cellcolor{gray!25} 0.105006    &\cellcolor{gray!25} 82.478540 	&\cellcolor{gray!25} 7.216383   \\
	\midrule
 	\texttt{SOCP$^M$} + \texttt{FBBT$^{NB}$}								 & \cellcolor{gray!25} 287	                      &\cellcolor{gray!25} 0.115036	  &	94.052187			            & 8.085631	\\
    \texttt{SOCP$^M$} + \texttt{0.5$\cdot$OV+0.5$\cdot$MP}                     & \cellcolor{gray!25} 289			              &	\cellcolor{gray!25} 0.104554	  &	90.193966			            & 7.653487	\\
	\texttt{FBBT$^{OB + NB}$} + \texttt{0.5$\cdot$OV+0.5$\cdot$MP}	         & \cellcolor{gray!25} 287 	                    &\cellcolor{gray!25} 0.103991	  &\cellcolor{gray!25} 81.727300	&\cellcolor{gray!25} 7.107938   \\
	\texttt{SOCP$^M$} + \texttt{FBBT$^{NB}$} + \texttt{0.5$\cdot$OV+0.5$\cdot$MP}	& \cellcolor{gray!25} 287		                &	\cellcolor{gray!25} 0.106330  &	94.880873			        &	7.883139	\\
    \bottomrule
    \end{tabular}}
\end{table}

\sisetup{
	round-mode = places,
	round-precision = 2,
	detect-weight=true,
	table-number-alignment = right,
	table-alignment-mode = none,
	table-text-alignment = right
	}

In Table \ref{tab:combined} we can see that, as expected from the previous results, the combinations that use approach \texttt{0.5$\cdot$OV+0.5$\cdot$MP} are the most competitive ones. Among them, probably \texttt{FBBT$^{OB + NB}$} + \texttt{0.5$\cdot$OV+0.5$\cdot$MP} is slightly superior to the rest, with the raw \texttt{0.5$\cdot$OV+0.5$\cdot$MP} performing very similarly. Interestingly, the combination \texttt{SOCP$^M$} + \texttt{0.5$\cdot$OV+0.5$\cdot$MP} is the one that solves most instances, despite not being so competitive in Time. \julio{It is worth noting that the approach that combines all enhancements, \texttt{SOCP$^M$} + \texttt{FBBT$^{NB}$} + \texttt{0.5$\cdot$OV+0.5$\cdot$MP}, is not the most competitive one, suggesting that there may be some interdependence between the different enhancements after all.}

The main takeaway message we have after the numerical results reported in this section, both for the individual enhancements and for their combinations, is that all enhancements seem to improve performance, although not to the same extent. Interestingly, the less sophisticated enhancement (and the easiest to implement), the strategy for the selection of the branching point, is the one with the (significantly) largest impact.

\section{Improving performance with machine learning}\label{sec:ml}
In this section we take the enhancements of Section~\ref{sec:domred} one step further. The goal is to study the potential for additional performance improvements by using machine learning techniques to predict the best combination of enhancements for a given instance. To this end, we also analyze the improvement on performance we would observe if we had an oracle signaling the best combination for each instance and, then, we use it as a reference to assess the quality of the learning. 


\subsection{Machine learning framework}
In \cite{Ghaddar:2023}, the authors present a machine learning framework designed to ``optimally'' choose between different configurations of a solver when confronting a new instance, given its underlying features. In order to do so, the machine learning model is first trained on a predefined set of instances, using both their features and the performance of the different configurations on each of them. \julio{Table~\ref{table:features} presents the full list of the input variables (features).\footnote{VIG and CMIG stand for two graphs that can be associated to any given polynomial optimization problem: \textit{variables intersection graph} and \textit{constraints-monomials intersection graph}, and whose precise definitions is given in \cite{Ghaddar:2023}.} They capture diverse characteristics of each instance and are a key ingredient of the machine learning framework.

\begin{table}[!htbp]
\centering
{\footnotesize
\renewcommand{\arraystretch}{0.6}
\setlength{\tabcolsep}{3pt}
    \begin{tabular}{ll}
      \toprule
      \multirow{5}{*}{Variables}   &  No. of variables, variance of the density of the variables \\
         & Average/median/variance of the ranges of the variables  \\
         & Average/variance of the no. of appearances of each variable\\
         & Pct. of variables not present in any monomial with degree greater than one \\
         & Pct. of variables not present in any monomial with degree greater than two \\
        \midrule
      Constraints &  No. of constraints, Pct. of equality/linear/quadratic constraints\\
    \midrule
    \multirow{3}{*}{Monomials}   &  No. of monomials\\
      & Pct. of linear/quadratic monomials, Pct. of linear/quadratic RLT variables\\
      & Average pct. of monomials in each constraint and in the objective function\\
      \midrule
      Coefficients   &  Average/variance of the coefficients\\
     \midrule
     \multirow{3}{*}{Other}   &  Degree and density of \ref{eq:PO}\\
     & No. of variables divided by no. of constraints/degree \\
     & No. of RLT variables/monomials divided by no. of constraints\\
     \midrule
    Graphs  & Density, modularity, treewidth, and transitivity of VIG and CMIG\\
    \bottomrule
    \end{tabular}}
\caption{Features used for the learning}
\label{table:features}
\end{table}
}

For our analysis we follow \cite{Ghaddar:2023}. We use Pace as the measure of performance (KPI) and use quantile regression forests for the learning \cite{Meinshausen:2006}, \julio{since in \cite{Ghaddar:2023} they proved to be slightly superior to  other techniques such as quantile generalized additive models and to stochastic gradient boosting for quantile regression}. \julio{We want to emphasize that we have deliberately employed a simple learning approach, without any customization or adaptation to our specific setting. The fact that this approach yields significant improvements, particularly in Section~\ref{sec:combinedlearning}, underscores the potential of such methodologies in this and related contexts. Yet, as we further discuss below, we acknowledge that it would be worth exploring how much farther these improvements might be pushed by using other learning methodologies such as deep learning and, in particular, graph neural networks \cite{gupta2020}, or enhancing the set of features with others that might be particularly relevant for domain reduction, \emph{e.g.} enriching the features related to ranges of the variables and their interactions with other features.}

One of the advantages of random forests, as well as other ensemble methods that use bagging, is that we the complete data set can be used to evaluate the model using the out-of-bag predictions, without randomly splitting the data into training and test sets.\footnote{As explained in \cite{Ghaddar:2023}, given a data set of size~$n$, ``a large number $m$ of bootstrap samples of size $n$ are obtained from the data set by sampling with replacement. Each bootstrap sample leaves out, on average, about 37\% of the observations (the left-out observations constitute the so-called out-of-bag sample). Then, $m$~individual trees are grown using these $m$ bootstrap samples. For each
observation in the data set, the corresponding out-of-bag prediction is obtained by taking into account only those trees fitted on bootstrap samples that leave out that particular observation''.} The learning is carried out on the same set of instances of the preceding section, using the statistical language R \cite{R-Core-Team:2023} with the library \emph{ranger} \cite{Wright:2017}. 

The reasons behind using Pace as KPI are the same ones that led to its introduction in \cite{Ghaddar:2023}: training can be performed on the whole set of instances at once. As discussed in Section~\ref{sec:testing}, Gap is not informative on instances solved by all approaches and Time is not informative on instances not solved by any approach. Pace, on the other hand, can discriminate between the performance of the different approaches regardless of the number of approaches that have solved each of the instances.

\julio{As mentioned above,} not only we follow the approach in \cite{Ghaddar:2023}, but we also use the same set of features and the same parameters for quantile regression. \julio{More precisely, borrowing from \cite{Ghaddar:2023}, ``learning is performed jointly on all sets of instances following the standard learning procedure. We randomly split the complete set of instances into two disjoint sets: the training set (70\%) and the test set (30\%). With the objective of obtaining a better performance, we gather the instances into ``families'' related to the groupings defined in the corresponding libraries, where those belonging to the same group share similar characteristics. The within-family proportion of instances is maintained through the splitting process. Moreover, in order to gain robustness, we construct 10 partitions of the dataset into training and test data and report aggregate results over all the partitions \ldots when performing quantile regression one has to specify the $\tau$-conditional quantile to be estimated.''\footnote{\julio{Given $\tau\in (0,1)$,  let $\hat{Q}_\tau(\mathbf{x})$ denote the estimated $\tau$-conditional quantile of $\NLBpace$ --~a normalization of $\LBpace$~-- for a given branching rule. This means that $(1-\tau)100\%$ of the time the performance of that rule is expected to be larger than $\hat{Q}_\tau(\mathbf{x})$.}} In the numerical analysis we set $\tau=0.3$, which proved to be the best performing value in \cite{Ghaddar:2023}, although the results were not particularly sensitive to this choice.}

In the results' tables in this section we report the performance of a special approach: \julio{\texttt{Oracle}}, which represents the performance obtained when we choose, for each instance, the best configuration according to pace, \emph{i.e.}, it's an ideal version representing what we would get if we had an oracle signaling the best approach for each instance. If, for a given instance, there are two options with the same pace, we choose the one with the smallest optimality gap.

\subsection{Impact of machine learning on the individual enhancements}

Table~\ref{tab:ML} summarizes the results of our learning methodology when applied individually to each one of the three different enhancements discussed in Section~\ref{sec:domred}. More precisely, for each enhancement, the goal is to predict which one of the considered approaches will have a better performance on a given instance. To simplify the tables of results, among the original approaches we just represent the best performing one for the given enhancement, \texttt{Best}. ``out-of-bag Q-RF'' is the result of choosing, for each instance, the approach selected by the machine learning framework. Finally, ``Improvement after learning'' is the percentage of relative improvement from \texttt{Best} to out-of-bag Q-RF and ``\julio{\texttt{Oracle}} improvement'' is the same but comparing \texttt{Best} to \julio{\texttt{Oracle}}.

\begin{table}[!htpb]
    \caption{Machine Learning impact on the different enhancements.}
    \label{tab:ML}
    \centering
    {\small \addtolength{\tabcolsep}{-2pt} 
    \begin{tabular}{l r S[round-precision=3] S S}
    \toprule
    \multicolumn{5}{c}{\fullref{sec:conic}}\\
    \midrule
	{\scriptsize (403 instances)} & {Solved {\scriptsize (403)}} & {Gap {\scriptsize (111)}} & {Time {\scriptsize (153)}} & {Pace {\scriptsize (268)}} \\
        \texttt{Best} (\texttt{Baseline}) & 286 & 0.132475 & 56.951878 & 6.560483 \\
        out-of-bag Q-RF & 286 & \cellcolor{gray!25} 0.130977 & \cellcolor{gray!25} 56.753052 & 6.567299 \\
        \julio{\texttt{Oracle}} & \cellcolor{gray!25} 287 & \cellcolor{gray!25} 0.127943 & \cellcolor{gray!25} 53.681877 & \cellcolor{gray!25} 6.191521 \\
        Improvement after learning & {$0.0\%$} &  \cellcolor{gray!25} {$-1.13\%$} &  \cellcolor{gray!25} {$-0.35\%$} & {$0.01\%$} \\
        \julio{\texttt{Oracle}} improvement & \cellcolor{gray!25} {$0.35\%$} & \cellcolor{gray!25} {$-3.42\%$} & \cellcolor{gray!25} {$-5.74\%$} & \cellcolor{gray!25} {$-5.62\%$} \\
    \midrule
    \multicolumn{5}{c}{\fullref{sec:bound-tightening}}\\
    \midrule
		{\scriptsize (403 instances)} & {Solved {\scriptsize (403)}} & {Gap {\scriptsize (111)}} & {Time {\scriptsize (134)}} & {Pace {\scriptsize (248)}} \\
        \texttt{Best} (\texttt{FBBT$^{OB+NB}$}) & 287 & 0.127502 & 98.320568 & 8.084024 \\
        out-of-bag Q-RF & \cellcolor{gray!25} 288 & \cellcolor{gray!25} 0.125127 & \cellcolor{gray!25} 94.717551 & \cellcolor{gray!25} 7.904663 \\
        \julio{\texttt{Oracle}} & \cellcolor{gray!25} 288 & \cellcolor{gray!25} 0.124292 & \cellcolor{gray!25} 89.653841 & \cellcolor{gray!25} 7.653371 \\
        Improvement after learning & \cellcolor{gray!25} {$0.35\%$} & \cellcolor{gray!25} {$-1.86\%$} & \cellcolor{gray!25} {$-3.66\%$} & \cellcolor{gray!25} {$-2.22\%$} \\
        \julio{\texttt{Oracle}} improvement & \cellcolor{gray!25} {$0.35\%$} & \cellcolor{gray!25} {$-2.52\%$} & \cellcolor{gray!25} {$-8.81\%$} & \cellcolor{gray!25} {$-5.33\%$} \\
    \midrule
    \multicolumn{5}{c}{\fullref{sec:branching-point}}\\
    \midrule
	{\scriptsize (403 instances)}	& {Solved {\scriptsize (403)}} & {Gap {\scriptsize (110)}} & {Time {\scriptsize (138)}} & {Pace {\scriptsize (250)}} \\
        \texttt{Best} (\texttt{0.5$\cdot$OV+0.5$\cdot$MP}) & 288 & 0.108060 & 96.746616 & 7.644139 \\
        out-of-bag Q-RF & 288 & 0.108497 & \cellcolor{gray!25} 94.853610 &  7.734806 \\
        \julio{\texttt{Oracle}} & \cellcolor{gray!25} 290 & \cellcolor{gray!25} 0.104173 & \cellcolor{gray!25} 82.658632 & \cellcolor{gray!25} 6.945399 \\
        Improvement after learning & {$0.0\%$} & {$0.4\%$} & \cellcolor{gray!25} {$-1.96\%$} & {$1.19\%$} \\
        \julio{\texttt{Oracle}} improvement & \cellcolor{gray!25} {$0.69\%$} & \cellcolor{gray!25} {$-3.6\%$} & \cellcolor{gray!25} {$-14.56\%$} & \cellcolor{gray!25} {$-9.14\%$} \\
    \bottomrule
    \end{tabular}}
\end{table}

The results in Table~\ref{tab:ML} are somewhat mixed. For the conic bound tightening from Section~\ref{sec:conic}, we see modest improvements in Gap and Time, but relatively far from the ones delivered by \julio{\texttt{Oracle}}. Learning seems to perform particularly well for the FBBT enhancements from Section~\ref{sec:bound-tightening}, with out-of-bag Q-RF significantly outperforming \texttt{Best} and getting around half-way of the improvements \julio{\texttt{Oracle}} would achieve. Finally, regarding the selection of the branching point from Section~\ref{sec:branching-point}, the learning does not seem to obtain much since, although it improves almost by 2\% in Time with respect to \texttt{Best}, its performance deteriorates in Gap and Pace by 0.4\% and 1.19\%, respectively. It seems as if the inherent randomness behind the impact of the branching point on the generated branch-and-bound trees is hard to learn upon. 

\subsection{Impact of machine learning on the combined enhancements}\label{sec:combinedlearning}

Despite the mixed results of the previous section, one would expect that, when adding all of the enhancements and their combinations together, there would be more potential for learning, given the additional richness and diversity of the underlying approaches/configurations. For the sake of exposition, learning is not carried out on all approaches and combinations, but just on those reported in Section~\ref{sec:combined}.

\begin{table}[!htpb]
    \caption{Machine Learning impact on all of the enhancements together.}
    \label{tab:MLcombined}
    \centering
    {\small \addtolength{\tabcolsep}{-2pt} 
    \begin{tabular}{l r S[round-precision=3] S S}
    \toprule
	{\scriptsize (403 instances)}	& {Solved {\scriptsize (403)}} & {Gap {\scriptsize (112)}} & {Time {\scriptsize (146)}} & {Pace {\scriptsize (257)}} \\
        \texttt{Best} (\texttt{FBBT$^{OB + NB}$} + \texttt{0.5$\cdot$OV+0.5$\cdot$MP}) & 287 & 0.103991 & 81.727300 & 7.107938 \\
        out-of-bag Q-RF & \cellcolor{gray!25} 290 & \cellcolor{gray!25} 0.103655 & \cellcolor{gray!25} 75.213011 & \cellcolor{gray!25} 6.885762 \\
        \julio{\texttt{Oracle}} & \cellcolor{gray!25} 291 & \cellcolor{gray!25} 0.099407 & \cellcolor{gray!25} 66.127683 & \cellcolor{gray!25} 6.275744 \\
        Improvement after learning & \cellcolor{gray!25} {$1.05\%$} & \cellcolor{gray!25} {$-0.32\%$} & \cellcolor{gray!25} {$-7.97\%$} & \cellcolor{gray!25} {$-3.13\%$} \\
        \julio{\texttt{Oracle}} improvement & \cellcolor{gray!25} {$1.39\%$} & \cellcolor{gray!25} {$-4.41\%$} & \cellcolor{gray!25} {$-19.09\%$} & \cellcolor{gray!25} {$-11.71\%$} \\
    \bottomrule
    \end{tabular}}
\end{table}

The results in Table~\ref{tab:MLcombined} show that machine learning performs fairly well. The learning approach outperforms \texttt{Best} (\texttt{FBBT$^{OB + NB}$} + \texttt{0.5$\cdot$OV+0.5$\cdot$MP}) according to all the metrics: it solves three more instances and obtains improvements of around 8\% and 3\% in Time and Gap. As expected, having a richer set of configurations to choose from also helps \julio{\texttt{Oracle}} to obtain larger improvements.

\begin{figure}[!htbp]
\caption{Percentages in which each version is selected by out-of-bag Q-RF and by \julio{\texttt{Oracle}}.}
\label{fig:comparison}
\centering
\includegraphics[scale=0.7, trim= 9 10 10 10, clip=true]{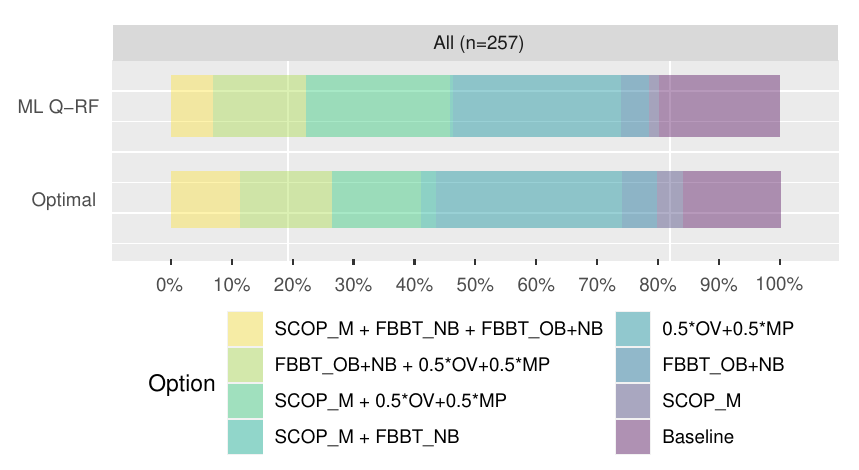}
\end{figure}

Figure~\ref{fig:comparison} represents, side by side, how often each of the eight configurations is selected by out-of-bag Q-RF and by \julio{\texttt{Oracle}}. We can see that out-of-bag Q-RF mimics quite well the behavior of \julio{\texttt{Oracle}}. Importantly, Figure~\ref{fig:comparison} also shows that a configuration with \texttt{SOCP$^M$} is chosen in around one third of the instances, demonstrating the potential of using conic-based relaxations in the OBBT stages of global optimization solvers. 

\begin{figure}[!htbp]
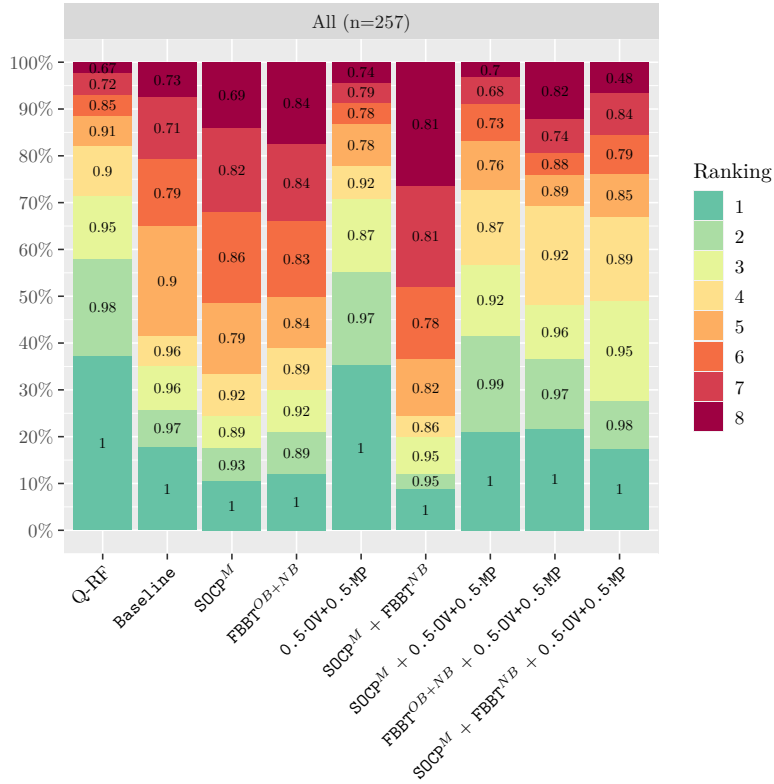

\caption{Domain reduction enhancements ranked from 1 (best) to 8 (worst), according to Pace.}
\label{fig:ranking}
\centering
\include{ranking_plot_combined.tex}
\vspace{-1cm}
\end{figure}

Despite the apparently nice behavior of out-of-bag Q-RF displayed in Figure~\ref{fig:comparison}, it might even be the case that it never chose the best configuration and just happened to use the different configurations with the right frequencies. Figure~\ref{fig:ranking} presents detailed information regarding the performance of the different configurations which, in particular, allows to assess how often each of them is the best one. The bar chart represents, for each configuration, the percentage of instances in which that configuration was the best one, the second best, and so on, for the eight ranking positions it might occupy depending on its performance relative to the rest.\footnote{Note that we have 8 configurations, in addition to out-of-bag Q-RF, and we just consider ranking positions from 1 to 8. This is because out-of-bag Q-RF always coincides with one of the eight configurations on which the learning is performed.} Moreover, the numbers inside the bars indicate how close that configuration is, on average, to the best one for the instances in which it occupies the corresponding ranking position. More specifically, these numbers go from zero to one and are computed by dividing the best (smallest) Pace among the different configurations by the Pace of the current one, and then taking the average over the instances in each ranking position. Thus, Figure~\ref{fig:ranking} allows to assess not only how often each approach is in each ranking position, but also how rapidly its performance deteriorates as it goes from the top positions to the bottom ones.

Figure~\ref{fig:ranking} shows that out-of-bag Q-RF comes on top more frequently than any other configuration, with \texttt{0.5$\cdot$OV+0.5$\cdot$MP} close behind. Yet, the main difference between them comes from the fact that the relative performance of out-of-bag Q-RF deteriorates more slowly as it moves from being the best configuration to the second or third best, for example. In Figure~\ref{fig:ranking} we can also see that the configurations that use \texttt{SOCP$^M$} are the ones in which performance tends to deteriorate more quickly, which is consistent with the discussion in Section~\ref{sec:conic} about their ``riskiness''. In this respect, since out-of-bag Q-RF uses \texttt{SOCP$^M$} in nearly one third of the instances, the fact that it performs so well even when it does not choose the top performing configuration suggests that it is successful at learning to avoid using \texttt{SOCP$^M$} in instances in which it does not perform well. Yet, the difference in performance between out-of-bag Q-RF and \julio{\texttt{Oracle}} also suggests that there is room for further improvement. A natural direction would be to enhance the set of features of the instances, with the goal of capturing additional aspects that might be relevant for domain reduction.

\section{Conclusions}\label{sec:conclusions}
We have studied the impact of various approaches to improve the effectiveness of domain reduction functionalities in the context of the RLT technique for solving polynomial optimization problems. Although all of these enhancements lead to performance improvements, the results showed that the overall impact may sometimes be relatively small. Interestingly, the approaches to select the branching point have the largest impact and they are, by far, the simplest and easiest to implement among all the enhancements under consideration. We believe that insights into the trade-offs between performance improvement and implementation costs are relevant to the design of more efficient optimization solvers, because they can guide the implementation efforts of state-of-the-art solver development teams toward the most promising directions. When a new technique or algorithmic enhancement is first studied, it is natural to analyze its impact in isolation to better quantify it. However, precisely because of this isolation, it is often difficult to determine which potential enhancements should be prioritized when designing or evolving a solver, as the reported impacts may not be easily comparable.

Following recent trends in the field, we built upon the learning framework developed in \cite{Ghaddar:2023} to demonstrate the potential of machine/statistical learning to enhance performance further. The promising results obtained in Section~\ref{sec:combinedlearning}, with the direct application of the methodology in  \cite{Ghaddar:2023}, suggest that one could achieve even better results with modifications tailored to specific applications. A relatively simple avenue for future research would be to extend the set of features used for learning to include additional features suitable for domain reduction.

\section*{Acknowledgments}
Project PID2021-124030NB-C32 funded by MICIU/AEI/10.13039/501100011033/ and by ERDF/EU. This research was also funded by Grupos de Referencia Competitiva ED431C-2021/24 from the Consellería de Cultura, Educación e Universidades, Xunta de Galicia. Ignacio Gómez-Casares acknowledges the support from the Spanish Ministry of Education through FPU grant 20/01555. 

\paragraph{Conflict of interest} The authors declare no competing interests.

\bibliography{references}
\bibliographystyle{ecta}

\appendix

\section{Additional results for FBBT enhancements}\label{app:FBBT}

Table~\ref{tab:zibbaronbtversions} presents a more comprehensive summary of the numerical results obtained when studying different variations of the FBBT approaches described in Section~\ref{sec:bound-tightening}. We study configurations in which the two approaches are applied at all nodes or just at the root node and at nodes of pre-specified depths. Note that, under the approach \texttt{FBBT$^{OB, 50 + NB, 10}$}, the cuts from FBBT$^{OB}$ are never applied by themselves, since multiples of~50 are also multiples of~10. In order to assess whether or not this might impact performance, we also studied \texttt{FBBT$^{OB, 50 + NB, 9}$}, under which both families of cuts are jointly applied at the same node much less often.

\julio{When looking at the performance of the approaches in which  \texttt{FBT$^{OB}$} or \texttt{FBT$^{NB}$} are applied at all nodes, we can see that the reductions on the sizes of the resulting trees do not compensate the computational overhead. In particular, when comparing \texttt{FBBT$^{OB, all}$} with \texttt{FBBT$^{OB, 25}$} and \texttt{FBBT$^{OB, 50}$}, the value
of Nodes shows that \texttt{FBBT$^{OB, all}$} yields slightly smaller trees than the other approaches but, on the other hand, \bttime becomes much larger. A similar observation applies to \texttt{FBBT$^{NB, all}$} with respect to \texttt{FBBT$^{NB, 5}$} and \texttt{FBBT$^{NB, 10}$}, although the increase in \bttime is less pronounced.}

\sisetup{
	round-mode = places,
	round-precision = 2,
	detect-weight=true,
	table-number-alignment = right,
	table-alignment-mode = none,
	table-text-alignment = right
	}

\begin{table}[!htpb]
    \caption{Performance of different versions of the FBBT enhanced with cuts from \cite{Gleixner:2017} and \cite{Ryoo:1996}.}
    \label{tab:zibbaronbtversions}
    \centering
    {\small \addtolength{\tabcolsep}{-2pt} 
    \begin{tabular}{l r S[round-precision=3] S S S | S}
    \toprule
        {\scriptsize (406 instances)} & {Solved {\scriptsize (403)}} & {Gap {\scriptsize (115)}} & {Time {\scriptsize (141)}} & {Pace {\scriptsize (255)}} & {Nodes {\scriptsize (281)}} & {\bttime {\scriptsize (403)}} \\
        {\texttt{Baseline}} & 286 & 0.111769 & 82.836185 & 7.167551 & 94.003800 & 280.991376 \\
        \midrule
        {\texttt{FBBT$^{OB, all}$}} & 281 & 0.155412 & 112.833067 & 9.835430 & 89.571381 & 819.167331 \\
        {\texttt{FBBT$^{OB, 25}$}} & 286 & 0.110492 & 84.867944 & 7.279889 & 91.891316 & 414.811738 \\
        {\texttt{FBBT$^{OB, 50}$}} & 287 & 0.106408 & 85.262235 & 7.194661 & 91.928627 & 300.358402 \\
        \midrule
        {\texttt{FBBT$^{NB, all}$}} & 287 & 0.109047 & 82.922871 & 7.130988 & 87.287363 & 329.931986 \\
        {\texttt{FBBT$^{NB, 5}$}} & 287 & 0.108293 & 83.252723 & 7.138033 & 88.299680 & 303.335822 \\
        {\texttt{FBBT$^{NB, 10}$}} & 287 & 0.107940 & 82.210871 & 7.087769 & 88.729864 & 299.392319 \\
        \midrule
        {\texttt{FBBT$^{OB, 25 + NB, 5}$}} & 284 & 0.116623 & 82.279138 & 7.152859 & 86.039915 & 429.977468 \\
        {\texttt{FBBT$^{OB, 25 + NB, 9}$}} & 284 & 0.116427 & 92.542588 & 7.631822 & 86.200926 & 429.160969 \\
        {\texttt{FBBT$^{OB, 25 + NB, 10}$}} & 284 & 0.115093 & 82.677573 & 7.171905 & 86.405732 & 428.969728 \\
        {\texttt{FBBT$^{OB, 50 + NB, 5}$}} & 286 & 0.109819 & 82.347894 & 7.068329 & 86.077355 & 321.100290 \\
        {\texttt{FBBT$^{OB, 50 + NB, 9}$}} & 286 & 0.110626 & 81.695155 & 7.034519 & 86.263637 & 317.244980 \\
        {\texttt{FBBT$^{OB, 50 + NB, 10}$}} & 287 & 0.107716 & 81.550860 & 7.029272 & 86.398599 & 317.244532 \\
    \bottomrule
    \end{tabular}}
\end{table}

\end{document}